\documentclass[a4paper]{amsart}

\usepackage[utf8]{inputenc}
\usepackage[english]{babel}

\usepackage[T1]{fontenc}
\usepackage{lmodern}  

\usepackage{my-math-symbol}
\usepackage{my-cleveref}
\usepackage{my-equation-numbering}
\usepackage{comment}
\usepackage{tikz}
\usetikzlibrary{cd}

\newcommand{\Cech}{\v{C}ech\xspace}
\newcommand{\reg}{\mathrm{reg}}
\newcommand{\sing}{\mathrm{sing}}

\newcommand{\rank}{\operatorname{rank}}
\newcommand{\Frac}{\operatorname{Frac}}
\newcommand{\Ext}{\operatorname{Ext}}
\newcommand{\calHom}{\calH om}

\usepackage[non-sorted-cites,initials,alphabetic,nobysame]{amsrefs}

\AtBeginDocument{%
	\def\MR#1{}
}

\title[Obstructions to pseudo-Einstein structures]{Complex singularities and obstructions to pseudo-Einstein structures: An expository note}
\author{Yuya Takeuchi}
\address{Division of Mathematics \\ Institute of Pure and Applied Sciences \\ University of Tsukuba
	\\ 1- 1- 1 Tennodai, Tsukuba, Ibaraki 305-8571 Japan}
\email{ytakeuchi@math.tsukuba.ac.jp, yuya.takeuchi.math@gmail.com}



\thanks{The author was supported by JSPS KAKENHI Grant Number JP25K17247.}

\begin{document}

\begin{abstract}
	This note explains the counterexamples to Lee's conjecture on pseudo-Einstein structures
	obtained in joint work with Yoshinori Gongyo, Kengo Hirachi, and Shuya Matsumoto.
	The examples are links of normal isolated singularities
	and occur in every odd real dimension at least three.
	Each link has vanishing first Betti number and vanishing integral first Chern class
	of its CR structure, yet admits no pseudo-Einstein structure.
	We explain how finite first homology relates the existence problem on the link
	to the canonical divisor at the singular point.
	For the higher-dimensional examples, we describe the key obstruction directly
	in terms of the first normal derivatives of the transition functions
	of a line bundle along an exceptional hypersurface.
\end{abstract}

\maketitle

\tableofcontents

\setcounter{section}{-1}
\section{Introduction}
\label{section:introduction}

Let $(M, T^{1, 0} M)$ be a closed strictly pseudoconvex CR manifold of CR dimension $n \geq 2$.
If $M$ admits a pseudo-Einstein structure,
the first Chern class of $T^{1, 0} M$ vanishes in $H^{2}(M, \bbR)$.
Lee~\cite{Lee1988} conjectured that this necessary condition was also sufficient.
Now let $M$ be a sufficiently small link of an isolated normal singularity $(X, o)$,
and write $X^{\ast} = X \setminus \{o\}$.
Then $M$ admits a pseudo-Einstein structure if and only if the restriction of the canonical bundle
$K_{X^{\ast}}$ to some inner collar of $M$ admits a flat Hermitian metric smooth up to $M$.
This characterization agrees with Lee's condition for $n \geq 2$.
Following Hirachi, Marugame, and Matsumoto~\cite{Hirachi-Marugame-Matsumoto2017}*{Proposition 2.6 and Definition 2.7},
we use it also in CR dimension one.

The following result was obtained in joint work with Gongyo, Hirachi, and Matsumoto.

\begin{theorem}[\cite{Gongyo-Hirachi-Matsumoto-Takeuchi2026}*{Theorem 1}]
\label{thm:main-counterexamples}
	For every integer $n \geq 1$, there is a closed connected strictly pseudoconvex CR manifold
	$(M, T^{1, 0} M)$ of real dimension $2 n + 1$ satisfying
	\begin{equation}
		b_{1}(M) = 0, \qquad c_{1}(T^{1, 0} M) = 0 \quad \text{in } H^{2}(M, \bbZ),
	\end{equation}
	and admitting no pseudo-Einstein structure.
	The manifold $M$ can be realized as a sufficiently small link
	of a normal isolated complex analytic singularity of complex dimension $n + 1$.
\end{theorem}

Independently, Dong and Wang~\cite{Dong-Wang2026}*{Theorem A}
constructed counterexamples in every odd real dimension at least five
using a related projective bundle construction.
The examples explained here have $b_{1}(M) = 0$ and also occur in CR dimension one.

This note explains the proof of \cref{thm:main-counterexamples}
for researchers in CR geometry who are unfamiliar with singular complex spaces.

The link between the CR geometry and the singularity is particularly direct when $H_{1}(M, \bbZ)$ is finite.
A flat Hermitian metric on the canonical bundle of an inner collar then has finite holonomy,
so some positive tensor power of that bundle has a parallel holomorphic frame.
By Hartogs extension, both the frame and its dual extend from the collar to the entire neighborhood bounded by $M$.
Their pairing remains equal to one, so the corresponding reflexive power of the canonical sheaf is locally free.
In divisor terminology, the canonical divisor is $\bbQ$-Cartier.
If $H_{1}(M, \bbZ) = 0$, the holonomy is already trivial, and the canonical divisor itself must be Cartier.
\Cref{section:complex-spaces,section:pseudo-einstein-canonical-divisors}
develop the analytic background needed for this argument.

The two constructions exploit these two conclusions.
In \cref{section:counterexample-integral-homology-sphere} we consider a normal surface singularity
of N\'emethi and Okuma~\cite{Nemethi-Okuma2017}*{Example 4.3}.
Its resolution graph shows that its link is an integral homology $3$-sphere,
while its non-Gorenstein property means that its canonical divisor is not Cartier.
The link therefore admits no pseudo-Einstein structure,
even though its CR structure $T^{1, 0} M$ is topologically trivial.

In higher dimensions,
the obstruction can be seen before the singularity is formed.
The construction produces a smooth hypersurface $D$ and a line bundle $F$ that is holomorphically trivial on $D$.
The first normal derivatives of its transition functions define a \Cech class on $D$.
This class is non-zero in our examples,
so no positive tensor power of $F$ is holomorphically trivial on any neighborhood of $D$.
\Cref{section:line-bundle-extensions-obstruction} develops this description of the first obstruction
directly from transition functions.
\Cref{section:projective-base,section:exceptional-hypersurface-projective-bundle}
construct the projective base and the projective bundle in which $D$ occurs.
\Cref{section:contraction-higher-cr-dimensions} then shows how contracting $D$ turns the obstruction
into the failure of the canonical divisor to be $\bbQ$-Cartier,
while the integral first Chern class on the link vanishes.

\section{Preliminaries on complex spaces}
\label{section:complex-spaces}

The purpose of this section is to explain the analytic language needed
when a complex manifold is replaced by a space with isolated singularities.
We emphasize the distinction between objects defined on the whole space
and objects defined only on its smooth part.
All complex spaces are assumed to be Hausdorff, second countable, and finite dimensional.
For general background, see Grauert--Remmert~\cite{Grauert-Remmert1984}*{Chapters 1, 4, and 6--8}
and \cite{Grauert-Remmert2004}*{Chapters A--B and IV--V}.

\subsection{Complex spaces and isolated singularities}
\label{subsection:complex-spaces-isolated-singularities}

The underlying set of a \emph{complex analytic space}, or simply a \emph{complex space},
is locally the common zero set of finitely many holomorphic functions.
The ideal generated by these functions determines the structure sheaf on this local model.
More precisely,
a local model consists of an open set $V \subset \bbC^{N}$,
a coherent ideal sheaf $\calI \subset \calO_{V}$,
and
\begin{equation}
	A = \Set{z \in V | f(z) = 0 \text{ for every local section } f \text{ of } \calI},
\end{equation}
equipped with a sheaf $\calO_{A} = (\calO_{V} / \calI)|_{A}$ of rings.
Here $\calO_{V}$ denotes the sheaf of holomorphic functions;
coherence is explained in \cref{subsection:coherent-sheaves-holomorphic-vector-bundles}.
The pair $(X, \calO_{X})$ is obtained by gluing such local models.
The sheaf $\calO_{X}$ is called the \emph{structure sheaf}.
Its stalk $\calO_{X, x}$ consists of germs of holomorphic functions at $x$.
Thus a complex space includes both a topological space and its local rings of functions.

A complex space is \emph{reduced} if no stalk contains a non-zero nilpotent element;
that is, $f^{m} = 0$ for some $m \geq 1$ implies $f = 0$.
Equivalently, the defining ideal in each local model is radical.
In this case,
it is the ideal of all holomorphic functions vanishing on $A$,
and sections of $\calO_{A}$ can be viewed as functions that locally extend holomorphically to $V$.
A continuous map between reduced complex spaces is holomorphic
if the pullback of every local holomorphic function is holomorphic.
In local embeddings,
this is the familiar condition that its coordinate functions are holomorphic.

A point $x \in X$ is \emph{regular}, or \emph{smooth},
if a neighborhood of $x$ is biholomorphic,
as a complex space,
to an open subset of some complex Euclidean space.
The regular locus is denoted by $X_{\reg}$;
its complement $X_{\sing}$ is the \emph{singular locus}.
If $X$ is reduced,
then the regular locus is dense and open,
and the singular locus is a closed analytic subset~\cite{Grauert-Remmert1984}*{6.2.2}.

The \emph{local complex dimension} $\dim_{x} X$ is the smallest integer $d \geq 0$
for which there are holomorphic functions $f_{1}, \dots, f_{d}$ near $x$
whose common zero set is $\{x\}$ in a sufficiently small neighborhood of $x$.
Equivalently,
\begin{equation}
	\dim_{x} X = \dim \calO_{X, x},
\end{equation}
where the right-hand side is the \emph{Krull dimension}:
the maximum length $r$ of a strictly increasing chain
$\mathfrak{p}_{0} \subsetneq \cdots \subsetneq \mathfrak{p}_{r}$ of prime ideals in $\calO_{X, x}$.
These definitions agree and are unchanged by replacing the local defining ideals by their radicals;
see \cite{Grauert-Remmert1984}*{5.1.1--5.1.3}.
At a regular point, this is the usual complex dimension of a manifold.
For a reduced space,
it is the largest manifold dimension occurring among regular points
in every sufficiently small neighborhood of $x$.
The \emph{dimension} of $X$ is
\begin{equation}
	\dim X \coloneqq \sup_{x \in X} \dim_{x} X.
\end{equation}
An irreducible component is a maximal analytic piece that cannot be expressed
as the union of two proper closed analytic subsets.
We say that $X$ has \emph{pure dimension} $d$
if all its irreducible components have complex dimension $d$.

An \emph{isolated singularity} at $o$ means that,
after shrinking $X$ around $o$,
\begin{equation}
	X_{\sing} = \{o\}.
\end{equation}
We write $X^{\ast} \coloneqq X \setminus \{o\}$ for this punctured neighborhood.
The notation $(X, o)$ denotes the germ of the space at $o$,
so assertions about $(X, o)$ are understood after possibly shrinking the neighborhood.

\subsection{Coherent sheaves and holomorphic vector bundles}
\label{subsection:coherent-sheaves-holomorphic-vector-bundles}

Let $(X, \calO_{X})$ be a complex space.
An $\calO_{X}$-module sheaf is a sheaf whose sections can be multiplied by holomorphic functions.
It is \emph{coherent} if it is locally generated by finitely many sections
and the relations among any finite collection of local sections are themselves locally finitely generated;
see \cite{Grauert-Remmert1984}*{A.3} for example.
The structure sheaf of a complex space is coherent by the coherence theorem;
see \cite{Grauert-Remmert1984}*{2.5.3}.
Consequently,
coherence of $\calF$ is equivalent to the existence,
near every point, of a finite presentation
\begin{equation}
	\calO_{U}^{\oplus a} \longrightarrow \calO_{U}^{\oplus b}
	\longrightarrow \calF|_{U} \longrightarrow 0.
\end{equation}
The first map is a matrix of holomorphic functions;
the presentation describes finitely many generators subject to finitely many relations.
Kernels, images, and cokernels of morphisms of coherent sheaves are coherent.
For these basic properties, see \cite{Grauert-Remmert1984}*{A.4.2}.

A holomorphic vector bundle $E$ of rank $r$ on $X$ is defined by local products
$U \times \bbC^{r}$ with holomorphic, invertible transition matrices.
Its sheaf $\calO(E)$ of holomorphic sections is \emph{locally free of rank $r$}:
locally it is isomorphic to $\calO_{X}^{\oplus r}$.
Conversely,
every locally free sheaf of finite rank arises in this way.
In particular, a holomorphic line bundle corresponds to a locally free sheaf of rank one,
also called an \emph{invertible sheaf}.

\subsection{Normality and holomorphic extension}
\label{subsection:normality-holomorphic-extension}

A complex space $X$ is \emph{normal}
if every local ring $\calO_{X, x}$ is an integrally closed integral domain.
An integral domain has no zero divisors;
being integrally closed means that an element of its fraction field
satisfying a monic polynomial with coefficients in the ring already belongs to the ring.
Normality is a local condition, stronger than reducedness.
Every complex manifold is normal.
A normal space is locally irreducible~\cite{Grauert-Remmert1984}*{6.4.2},
and its singular locus has complex codimension at least two~\cite{Grauert-Remmert1984}*{6.5.3}.

There is an analytic interpretation of normality.
A function on $X_{\reg}$ is called \emph{weakly holomorphic}
if it is holomorphic there and locally bounded near every point of $X$.
A reduced complex space is normal if and only if, for every open subset $U \subset X$,
every weakly holomorphic function on $U_{\reg}$
extends to an element of $\calO_{X}(U)$~\cite{Grauert-Remmert1984}*{7.4.2}.
Thus normality says that locally bounded holomorphic functions on the smooth locus
are already genuine holomorphic functions on the space.

Normality also gives extension without a boundedness assumption across smaller analytic sets.

\begin{theorem}[Riemann extension theorem~\cite{Grauert-Remmert1984}*{7.4.2}]
\label{thm:normal-extension}
	Let $X$ be a normal complex space
	and let $A \subset X$ be a closed analytic subset of complex codimension at least two.
	If $j \colon X \setminus A \hookrightarrow X$ is the inclusion,
	then
	\begin{equation}
		\calO_{X} \cong j_{\ast} \calO_{X \setminus A}.
	\end{equation}
	Equivalently,
	for every open set $U \subset X$,
	restriction induces an isomorphism
	\begin{equation}
		\calO_{X}(U) \stackrel{\cong}{\longrightarrow} \calO_{X}(U \setminus A).
	\end{equation}
\end{theorem}

Let $\calF$ be a coherent sheaf on a normal complex space $X$.
Then the local ring $R_{x} \coloneqq \calO_{X, x}$ is an integral domain.
Let $Q_{x} \coloneqq \Frac R_{x}$ be its fraction field.
We define
\begin{equation}
	\rank_{x} \calF
	\coloneqq \dim_{Q_{x}} \rbra*{\calF_{x} \otimes_{R_{x}} Q_{x}}.
\end{equation}

\begin{lemma}
	The rank $\rank_{x} \calF$ of a coherent sheaf $\calF$ is locally constant.
\end{lemma}

\begin{proof}
	Take a finite presentation
	\begin{equation}
		\calO_{U}^{\oplus a} \stackrel{B}{\longrightarrow} \calO_{U}^{\oplus b}
		\longrightarrow \calF|_{U} \longrightarrow 0
	\end{equation}
	on a sufficiently small connected neighborhood $U$ of $x$.
	Tensoring the stalk sequence with $Q_{x}$ gives
	\begin{equation}
		\rank_{x} \calF = b - \rank_{Q_{x}} B_{x},
	\end{equation}
	where $B_{x} \colon Q_{x}^{\oplus a} \to Q_{x}^{\oplus b}$ is the $Q_{x}$-linear map induced by $B$.
	Since the rank of a linear map is determined by its minors,
	the identity theorem implies $\rank_{Q_{x}} B_{x}$ is locally constant,
	and so is $\rank_{x} \calF$.
\end{proof}

We say that $\calF$ has rank $r$
if $\rank_{x} \calF = r$ for every $x \in X$.
Equivalently,
$\calF$ is locally free of rank $r$ on a dense open subset of $X$.

An element $s \in \calF_{x}$ is a \emph{torsion element}
if $a s = 0$ for some non-zero $a \in R_{x}$.
The sheaf $\calF$ is \emph{torsion-free}
if the only torsion element of each stalk $\calF_{x}$ is zero.
Equivalently, the natural map
\begin{equation}
	\calF_{x} \longrightarrow \calF_{x} \otimes_{R_{x}} Q_{x}
\end{equation}
is injective for every $x \in X$.
In particular,
a torsion-free coherent sheaf of rank zero is the zero sheaf,
since each stalk embeds in a zero-dimensional $Q_{x}$-vector space.

\subsection{Reflexive sheaves and \texorpdfstring{$\bbQ$}{Q}-Cartier canonical divisors}
\label{subsection:reflexive-sheaves-canonical-divisors}

Throughout this subsection,
$X$ is normal and of pure dimension $d$.
For a coherent sheaf $\calF$,
set
\begin{equation}
	\calF^{\vee} \coloneqq \calHom_{\calO_{X}}(\calF,\calO_{X}),
	\qquad
	\calF^{\vee \vee} \coloneqq (\calF^{\vee})^{\vee}.
\end{equation}
Note that both $\calF^{\vee}$ and $\calF^{\vee \vee}$ are coherent~\cite{Grauert-Remmert1984}*{A.4.4}.
The sheaf $\calF$ is \emph{reflexive}
if the natural evaluation map $\calF \to \calF^{\vee\vee}$ is an isomorphism.
Every locally free sheaf is reflexive,
but a reflexive sheaf need not be locally free.
A reflexive sheaf is torsion-free;
dualizing a local finite presentation of $\calF^{\vee}$ embeds $\calF$ in a finite free sheaf,
which is torsion-free.

\begin{lemma}
	Let $\calF$ be a coherent sheaf on a normal complex space $X$.
	Then $\calF^{\vee}$ is reflexive.
\end{lemma}

\begin{proof}
	Work at a stalk and write $R = \calO_{X, x}$ and $M = \calF_{x}$,
	with $M^{\vee} = \Hom_{R}(M, R)$.
	If $\delta_{M} \colon M \to M^{\vee \vee}$ is the evaluation map,
	then
	\begin{equation}
		\delta_{M}^{\vee} \circ \delta_{M^{\vee}} = \id_{M^{\vee}}.
	\end{equation}
	Thus $M^{\vee \vee \vee} \cong M^{\vee} \oplus \Ker \delta_{M}^{\vee}$.
	Since $M$ is finitely presented, dualization commutes with tensoring with $Q = \Frac R$:
	\begin{equation}
		\Hom_{R}(M, R) \otimes_{R} Q
		\cong \Hom_{Q}(M \otimes_{R} Q, Q).
	\end{equation}
	It follows that $M^{\vee}$ and $M^{\vee \vee \vee}$ have the same rank.
	The kernel has rank zero and is torsion-free as a submodule of a dual module over a domain,
	so it must vanish.
	Hence $\delta_{M^{\vee}}$ is an isomorphism.
\end{proof}

Since $\calF^{\vee}$ is coherent,
the preceding lemma implies that $\calF^{\vee \vee}$ is also reflexive.
The sheaf $\calF^{\vee \vee}$ is called the \emph{reflexive hull} of $\calF$.

A reflexive sheaf is determined by its restriction outside any analytic subset of codimension at least two:

\begin{lemma}
	Let $\calF$ be a coherent reflexive sheaf on $X$
	and let $A \subset X$ be a closed analytic subset of complex codimension at least two.
	If $j \colon X \setminus A \hookrightarrow X$ is the inclusion,
	then
	\begin{equation}
	\label{eq:reflexive-extension}
		\calF \cong j_{\ast}(\calF|_{X \setminus A}).
	\end{equation}
	Equivalently,
	for every open set $U \subset X$,
	restriction induces an isomorphism
	\begin{equation}
		\calF(U) \stackrel{\cong}{\longrightarrow} \calF(U \setminus A).
	\end{equation}
\end{lemma}

\begin{proof}
	For a sufficiently small neighborhood $U$ of each point,
	take a finite presentation
	\begin{equation}
		\calO_{U}^{\oplus a} \longrightarrow \calO_{U}^{\oplus b}
		\longrightarrow \calF^{\vee}|_{U} \longrightarrow 0
	\end{equation}
	of the coherent sheaf $\calF^{\vee}$.
	Since $\calHom_{\calO_{U}}(\cdot, \calO_{U})$ is left-exact,
	we have
	\begin{equation}
		0 \longrightarrow \calF|_{U} \longrightarrow \calO_{U}^{\oplus b}
		\stackrel{B}{\longrightarrow} \calO_{U}^{\oplus a}.
	\end{equation}
	Thus a section of $\calF$ is a $b$-tuple $(f_{1}, \dots , f_{b})$
	of holomorphic functions satisfying $B(f_{1}, \dots , f_{b}) = 0$.
	Take a section $s = (f_{1}, \dots , f_{b}) \in \calF(U \setminus A)$.
	Extend each component $f_{i}$ to $\wtf_{i} \in \calO_{X}(U)$ by \cref{thm:normal-extension}.
	The extended $b$-tuple $(\wtf_{1}, \dots , \wtf_{b})$ still satisfies
	the equations $B(\wtf_{1}, \dots , \wtf_{b}) = 0$ by the identity theorem,
	and hence defines a section of $\calF$.
	Uniqueness also follows from the identity theorem.
\end{proof}

The cotangent sheaf $\Omega_{X}^{1}$ is the sheaf generated locally
by differentials of holomorphic functions,
subject to linearity, the Leibniz rule, and the holomorphic chain rule.

\begin{lemma}
	The cotangent sheaf $\Omega_{X}^{1}$ is coherent.
\end{lemma}

\begin{proof}
	Fix $x \in X$ and choose a neighborhood $U$ of $x$
	and a closed holomorphic embedding
	\begin{equation}
		\iota \colon U \longrightarrow V \subset \bbC^{N}.
	\end{equation}
	Let $\calI \subset \calO_{V}$ be the defining ideal sheaf.
	After shrinking $V$,
	choose holomorphic generators
	\begin{equation}
		\calI = (f_{1}, \dots, f_{s}).
	\end{equation}
	Such a finite set of generators exists
	because $\calI$ is coherent~\cite{Grauert-Remmert1984}*{4.2}.
	The local definition of the cotangent sheaf is
	\begin{equation}
		\Omega_{U}^{1}
		= \rbra*{\Omega_{V}^{1} / (\calI \Omega_{V}^{1} + \calO_{V} \cdot d \calI)}|_{U}.
	\end{equation}
	Since
	\begin{equation}
		\Omega_{V}^{1} = \bigoplus_{j = 1}^{N} \calO_{V} d z_{j},
	\end{equation}
	we have
	\begin{equation}
		\Omega_{U}^{1}
		= \rbra*{\bigoplus_{j = 1}^{N} \calO_{U} d z_{j}} \Bigg/ \sum_{a = 1}^{s} \calO_{U}
			\rbra*{\sum_{j = 1}^{N} \left. \pdv{f_{a}}{z_{j}} \right|_{U} d z_{j}}.
	\end{equation}
	This gives the finite presentation
	\begin{equation}
		\calO_{U}^{\oplus s} \longrightarrow \calO_{U}^{\oplus N}
		\longrightarrow \Omega_{U}^{1} \longrightarrow 0.
	\end{equation}
	In particular,
	$\Omega_{U}^{1}$ is the cokernel of a morphism between coherent sheaves,
	and so it is also coherent.
\end{proof}

On $X_{\reg}$ it is the usual holomorphic cotangent bundle,
but on $X$ it may have torsion and need not be locally free.
We write $\Omega_{X}^{p} = \bigwedge^{p} \Omega_{X}^{1}$,
which is coherent by \cite{Grauert-Remmert2004}*{A.2.5},
and
\begin{equation}
	\Omega_{X}^{[p]}
	\coloneqq (\Omega_{X}^{p})^{\vee \vee}
	\cong j_{\ast} \Omega_{X_{\reg}}^{p},
\end{equation}
where $j \colon X_{\reg} \hookrightarrow X$ is the inclusion.
These are the sheaves of \emph{reflexive differential forms}.
In particular,
the \emph{reflexive canonical sheaf} is
\begin{equation}
\label{eq:canonical-sheaf}
	K_{X}
	\coloneqq \Omega_{X}^{[d]}
	\cong j_{\ast} \calO(K_{X_{\reg}}).
\end{equation}
It is a coherent reflexive sheaf of rank one.
Its sections are holomorphic $d$-forms on the smooth locus,
interpreted as sections of a sheaf on the whole space.

For the purposes of this note,
the relevant divisor condition concerns only the canonical sheaf.
For a positive integer $m$,
its \emph{reflexive $m$-th power} is
\begin{equation}
	K_{X}^{[m]}
	\coloneqq (K_{X}^{\otimes m})^{\vee \vee}
	\cong j_{\ast} \calO(K_{X_{\reg}}^{\otimes m}).
\end{equation}
We say that the \emph{canonical divisor is $\bbQ$-Cartier at $x$}
if $K_{X}^{[m]}$ is invertible on a neighborhood of $x$ for some positive integer $m$.

Locally,
we also use $K_{X}$ to denote a canonical Weil divisor.
A non-zero local meromorphic top-degree form on the smooth locus determines such a divisor
by its orders of zeros and poles along the codimension-one analytic subsets.
The associated rank-one reflexive sheaf is $\calO_{X}(K_{X}) \cong K_{X}$.
The statement that $m K_{X}$ is Cartier means precisely that
$\calO_{X}(m K_{X}) \cong K_{X}^{[m]}$ is a line bundle.
Thus the condition can be formulated entirely in terms of the reflexive canonical sheaf,
without choosing a canonical divisor or developing the theory of arbitrary divisors.
It asserts local freeness of a reflexive power, not global triviality
or vanishing of its first Chern class.

Let $(X, o)$ be a normal isolated singularity of dimension at least two.
There is a useful equivalent formulation after shrinking around $o$:
the canonical divisor $K_{X}$ is $\bbQ$-Cartier at $o$ if and only if
$K_{X^{\ast}}^{\otimes m}$ is holomorphically trivial for some $m \geq 1$.
Indeed,
an invertible $K_{X}^{[m]}$ has a local frame near $o$,
whose restriction trivializes the bundle on $X^{\ast}$.
Conversely,
a trivialization on $X^{\ast}$ gives
\begin{equation}
	K_{X}^{[m]}
	\cong j_{\ast} \calO_{X^{\ast}}
	\cong \calO_{X}
\end{equation}
by \cref{thm:normal-extension}.
This criterion concerns a sufficiently small punctured neighborhood;
it does not assert triviality on the smooth locus of an arbitrary global space.

\subsection{Stein spaces and Hartogs extension}
\label{subsection:stein-spaces-hartogs-extension}

Stein spaces play the role of domains with sufficiently many global holomorphic functions.
A reduced complex space $X$ is \emph{Stein}
if it is holomorphically separable and holomorphically convex.
The first condition means that distinct points can be separated by a global holomorphic function.
The second means that for every compact set $K \subset X$, its holomorphically convex hull
\begin{equation}
	\whK_{X}
	\coloneqq \Set{x \in X | \abs{f(x)} \leq \sup_{K} \abs{f} \text{ for every } f \in \calO_{X}(X)}
\end{equation}
is compact.
Closed analytic subspaces of Stein spaces are Stein;
in particular,
closed analytic subsets of $\bbC^{N}$ are Stein.
Every point of a reduced complex space has a basis of Stein neighborhoods,
obtained by intersecting the image of a local closed embedding with sufficiently small balls or polydiscs.

\begin{theorem}[Cartan's theorems A and B~\cite{Grauert-Remmert2004}*{IV.1.1 and V.4.2}]
\label{thm:cartan}
	Let $X$ be a Stein space and let $\calF$ be a coherent sheaf on $X$.
	Then every stalk $\calF_{x}$ is generated over $\calO_{X, x}$ by germs of global sections,
	and $H^{q}(X, \calF) = 0$ for every $q \geq 1$.
\end{theorem}

The following global Hartogs theorem concerns compact holes,
which need not be analytic subsets.

\begin{theorem}[\cite{Kasahara1965}*{Theorem 3}]
\label{thm:stein-hartogs}
	Let $W_{0}$ be a connected normal Stein space of dimension at least two.
	Let $W \subset W_{0}$ be a connected open set
	and let $K \subset W$ be compact.
	If $W \setminus K$ is connected,
	every holomorphic function on $W \setminus K$
	extends uniquely to a holomorphic function on $W$.
	Similarly,
	every meromorphic function on $W \setminus K$
	extends uniquely to a meromorphic function on $W$.
\end{theorem}

As an application of this theorem,
we can extend sections of coherent reflexive sheaves across compact holes.

\begin{lemma}[Hartogs extension for reflexive sheaves]
\label{lem:reflexive-compact-hartogs}
	Let $W_{0}$ be a connected normal Stein space of dimension at least two,
	let $W \Subset W_{0}$ be a relatively compact connected open set,
	and let $K \subset W$ be compact
	such that $C \coloneqq W \setminus K$ is non-empty and connected.
	If $\calE$ is a coherent reflexive sheaf on $W_{0}$,
	then restriction induces an isomorphism
	\begin{equation}
		\calE(W) \stackrel{\cong}{\longrightarrow} \calE(C).
	\end{equation}
\end{lemma}

\begin{proof}
	Since $\ovW$ is compact and $\calE^{\vee}$ is coherent,
	we can choose finitely many global sections
	which generate $\calE^{\vee}$ on a neighborhood of $\ovW$.
	They define a morphism $\calO_{W_{0}}^{\oplus b} \to \calE^{\vee}$ on all of $W_{0}$.
	Its kernel is coherent,
	so a second application of Cartan's theorem A gives a presentation
	\begin{equation}
		\calO_{W_{0}}^{\oplus a} \longrightarrow \calO_{W_{0}}^{\oplus b}
		\longrightarrow \calE^{\vee} \longrightarrow 0,
	\end{equation}
	which is exact on a neighborhood of $\ovW$.
	Since $\calHom_{\calO_{W}}(\cdot, \calO_{W})$ is left-exact,
	we have
	\begin{equation}
		0 \longrightarrow \calE \longrightarrow \calO_{W}^{\oplus b}
		\stackrel{B}{\longrightarrow} \calO_{W}^{\oplus a}.
	\end{equation}
	Thus a section of $\calE$ is a $b$-tuple $(f_{1}, \dots , f_{b})$
	of holomorphic functions satisfying $B(f_{1}, \dots , f_{b}) = 0$.
	Take a section $s = (f_{1}, \dots , f_{b}) \in \calE(C)$.
	Extend each component $f_{i}$ to $\wtf_{i} \in \calO_{W}(W)$ by \cref{thm:stein-hartogs}.
	The extended $b$-tuple $(\wtf_{1}, \dots , \wtf_{b})$ still satisfies
	the equations $B(\wtf_{1}, \dots , \wtf_{b}) = 0$ by the identity theorem,
	and hence defines a section of $\calE$.
	Uniqueness also follows from the identity theorem.
\end{proof}

\subsection{Links of singularities and contraction of exceptional sets}
\label{subsection:links-contractions}

Let $(X, o)$ be a normal isolated singularity of dimension $n + 1 \geq 2$,
and choose a local embedding in $\bbC^{N}$ sending $o$ to the origin.
Denote by $\bbB_{r}$ an open ball of radius $r$ centered at $0$ in $\bbC^{N}$,
and set $\bbS_{r} \coloneqq \bdry \bbB_{r}$.
On a sufficiently small representative,
consider the function
\begin{equation}
	d_{X} \colon X \longrightarrow \bbR; \qquad z \longmapsto \norm{z}^{2}.
\end{equation}
By the local conical structure theorem \cite{Seade2019}*{Theorem 3.1},
there exists $\varepsilon > 0$ such that
$\bbS_{r}$ intersects $X^{\ast}$ transversely for any $0 < r \leq \varepsilon$.
This intersection
\begin{equation}
	M_{r} \coloneqq X^{\ast} \cap \bbS_{r}
\end{equation}
is a closed smooth manifold of real dimension $2 n + 1$,
called the \emph{link} of the singularity.
Moreover,
there exists a smooth one-parameter family of diffeomorphisms
\begin{equation}
	\gamma_{t} \colon \bbS_{\varepsilon} \longrightarrow \bbS_{\varepsilon - t}, \qquad
	t \in \clop{0}{\varepsilon}
\end{equation}
such that $\gamma_{0}$ is the identity and $\gamma_{t}(M_{\varepsilon}) = M_{\varepsilon - t}$.
Furthermore,
the pair $(\bbB_{\varepsilon}, X \cap \bbB_{\varepsilon})$ is homeomorphic to
the open cone over $(\bbS_{\varepsilon}, M_{\varepsilon})$.
In particular,
$X \cap \bbB_{\varepsilon}$ is connected.
The normality of $X$ implies that
restriction induces an isomorphism
\begin{equation}
	\calO_{X}(X \cap \bbB_{\varepsilon}) \stackrel{\cong}{\longrightarrow}
	\calO_{X}(X^{\ast} \cap \bbB_{\varepsilon}).
\end{equation}
If $X^{\ast} \cap \bbB_{\varepsilon}$ were disconnected,
it would admit a non-constant holomorphic function $f$ with values in $\{0, 1\}$.
Its holomorphic extension $\wtf$ to $X \cap \bbB_{\varepsilon}$ would satisfy
$\wtf^{2} = \wtf$ by continuity,
and hence would be constant since $X \cap \bbB_{\varepsilon}$ is connected.
This is a contradiction,
so $X^{\ast} \cap \bbB_{\varepsilon}$ is connected.
The conical structure identifies $X^{\ast} \cap \bbB_{\varepsilon}$
with $M_{\varepsilon} \times \opop{0}{\varepsilon}$,
with the second coordinate given by the radius.
It follows that all the links $M_{r}$ are connected.
This radial product structure also shows that,
for $0 < r^{\prime} < r \leq \varepsilon$,
the inclusion
\begin{equation}
	X^{\ast} \cap \bbB_{r^{\prime}} \hookrightarrow X^{\ast} \cap \bbB_{r}
\end{equation}
is a homotopy equivalence.

The link carries more structure than its underlying smooth manifold.
As a real hypersurface in $X^{\ast}$,
it has the induced CR structure
\begin{equation}
	T^{1, 0} M_{r}
	= T^{1, 0} X^{\ast}|_{M_{r}}
	\cap (T M_{r} \otimes \bbC).
\end{equation}
The function $d_{X}$ is strictly plurisubharmonic on $X^{\ast}$,
so its regular level hypersurfaces are strictly pseudoconvex.
Thus the link is naturally a strictly pseudoconvex CR manifold.

An important construction of Stein spaces with normal isolated singularities
is the Remmert reduction of a strictly pseudoconvex domain.
Let $\Omega$ be a relatively compact domain with smooth strictly pseudoconvex boundary
in a complex manifold.
Then there is a proper holomorphic map $\pi \colon \Omega \to W$ onto a normal Stein space such that
\begin{equation}
	\pi_{\ast} \calO_{\Omega} = \calO_{W}.
\end{equation}
The union $E$ of the positive-dimensional irreducible compact analytic subsets of $\Omega$
is a compact analytic set with finitely many connected components.
The map contracts each of these components to a point
and induces a biholomorphism $\Omega \setminus E \cong W \setminus \pi(E)$.
In particular,
$W$ is smooth away from finitely many points.
For this reduction,
see \cite{Grauert1962}*{S\"{a}tze 1, 4--5 in Section 2}.

In the construction below, we start with a compact smooth hypersurface in a complex manifold
and seek a map that contracts it to a point.
A holomorphic line bundle is called \emph{negative}
if it admits a smooth Hermitian metric whose curvature form is negative definite.
Equivalently, its dual admits a Hermitian metric with positive definite curvature.
In particular, the dual of an ample line bundle on a smooth projective manifold is negative.
The following theorem of Grauert gives a sufficient condition for such a contraction
in terms of the negativity of the normal bundle of the hypersurface.

\begin{theorem}[\cite{Grauert1962}*{Satz 8 in Section 3}]
\label{thm:grauert-contraction}
	Let $D$ be a compact connected smooth hypersurface of positive dimension in a complex manifold $Y$.
	If its normal line bundle $N_{D / Y}$ is negative,
	then $D$ is an exceptional analytic set.
	In particular,
	there exist a neighborhood $\Omega$ of $D$
	and a proper holomorphic map
	\begin{equation}
		q \colon \Omega \longrightarrow W
	\end{equation}
	onto a normal Stein space,
	with $q^{- 1}(o) = D$ for a point $o \in W$,
	such that $\Omega \setminus D \to W \setminus \{o\}$ is biholomorphic.
\end{theorem}

\section{Pseudo-Einstein structures and canonical divisors}
\label{section:pseudo-einstein-canonical-divisors}

\subsection{Pseudo-Einstein structures and flat Hermitian metrics}
\label{subsection:pseudo-einstein-structures-flat-hermitian-metrics}

Let $(M, T^{1, 0} M)$ be a closed strictly pseudoconvex CR manifold of dimension $2 n + 1 \geq 3$.
A contact form $\theta$ is called a \emph{pseudo-Einstein structure} if
\begin{equation}
	\begin{cases}
		\Ric_{\alpha \ovxb} = \dfrac{\Scal}{n} h_{\alpha \ovxb},
			& n \geq 2, \\
		\nabla_{1} \Scal - i \nabla^{1} A_{1 1} = 0,
			& n = 1.
	\end{cases}
\end{equation}
For $n \geq 2$, this is Lee's original condition~\cite{Lee1988}.
Hirachi, Marugame, and Matsumoto gave a characterization of the pseudo-Einstein condition
in terms of the CR canonical bundle.
The CR canonical bundle is
\begin{equation}
	K_{M}
	\coloneqq \Set{\zeta \in \bigwedge\nolimits^{\! n + 1}(T^{\ast} M \otimes \bbC) |
	Z \contr \zeta = 0 \text{ for all } Z \in T^{0, 1} M}.
\end{equation}
A contact form $\theta$ defines a Hermitian metric $h_{\theta}$ on $K_{M}$,
and conversely,
a Hermitian metric on $K_{M}$ yields a contact form on $M$.
The CR canonical bundle with a Hermitian metric has a canonical metric-preserving connection $D$,
called the Tanaka connection.
The condition that $\theta$ is pseudo-Einstein is equivalent to flatness of $D$%
~\cite{Hirachi-Marugame-Matsumoto2017}*{Proposition 2.6}.

Assume that $M$ is a real hypersurface in an $(n + 1)$-dimensional complex manifold $N$.
In this case,
pullback of ambient forms identifies $K_{M}$ with $K_{N}|_{M}$.
Under this identification,
$\theta$ is pseudo-Einstein if and only if
$h_{\theta}$ admits a flat extension to a collar on the pseudoconvex side, smooth up to $M$%
~\cite{Hirachi-Marugame-Matsumoto2017}*{Proposition 2.6 and Definition 2.7}.

Let $\rho$ be a smooth real-valued defining function for $M$ in $N$.
Then $\partial \rho$ induces the exact sequence of smooth complex vector bundles
\begin{equation}
	0 \longrightarrow T^{1, 0} M \longrightarrow T^{1, 0} N|_{M}
	\stackrel{\partial \rho}{\longrightarrow} M \times \bbC \longrightarrow 0.
\end{equation}
Taking determinants and then first Chern classes gives
\begin{equation}
	c_{1}(K_{N}|_{M}) = - c_{1}(T^{1, 0} M).
\end{equation}
In particular,
if $M$ admits a pseudo-Einstein structure,
then $c_{1}(T^{1, 0} M) = 0$ in $H^{2}(M, \bbR)$.
For $n \geq 2$, Lee's conjecture asks whether the converse holds.

\subsection{Canonical divisors and finite holonomy}
\label{subsection:canonical-divisors-finite-holonomy}

Let $M$ be a sufficiently small link of a normal isolated singularity
$(X, o)$ of complex dimension $n + 1 \geq 2$.
For links,
we apply the preceding criterion in the smooth locus $N = X^{\ast}$.
We first show that
the $\bbQ$-Cartierness of $K_{X}$ at $o$ implies
the existence of pseudo-Einstein structures on links;
here we do not need additional topological assumptions.

\begin{lemma}
	Let $(X, o)$ be a normal isolated singularity of complex dimension $n + 1 \geq 2$
	whose canonical divisor $K_{X}$ is $\bbQ$-Cartier at $o$.
	Then any sufficiently small link of $(X, o)$ admits a pseudo-Einstein structure.
\end{lemma}

\begin{proof}
	Assume that $K_{X}^{[m]}$ is invertible at $o$.
	Choose a local frame $s$ on a sufficiently small representative,
	and restrict it to the smooth punctured neighborhood $X^{\ast}$.
	For a local holomorphic frame $e$ of $K_{X^{\ast}}$,
	write $s = \varphi e^{m}$ with $\varphi$ holomorphic and nowhere-vanishing.
	The requirement that $s$ have unit norm defines a Hermitian metric $h$ on $K_{X^{\ast}}$:
	\begin{equation}
		\norm{e}_{h}^{2} = \abs{\varphi}^{- 2 / m}.
	\end{equation}
	The metric $h$ is independent of the choice of $e$,
	and its curvature vanishes since
	\begin{equation}
		\partial \delb \log \norm{e}_{h}^{2}
		= - \frac{1}{m} \partial \delb \log \abs{\varphi}^{2}
		= 0.
	\end{equation}
	Choosing the link inside this representative ensures that
	the metric is defined and smooth on a neighborhood of the link.
	It therefore gives a pseudo-Einstein structure.
\end{proof}

If the first Betti number of the link is equal to zero,
then the converse also holds.

\begin{proposition}
\label{prop:finite-holonomy-canonical}
	Let $M$ be a sufficiently small link of
	a normal isolated singularity $(X, o)$ of dimension at least two.
	If $H_{1}(M, \bbZ)$ is finite and $K_{X}$ is not $\bbQ$-Cartier at $o$,
	then $M$ admits no pseudo-Einstein structure.
	If $H_{1}(M, \bbZ) = 0$ and $K_{X}$ is not Cartier at $o$,
	the same conclusion holds.
\end{proposition}

\begin{proof}
	Fix a local embedding in $\bbC^{N}$ sending $o$ to the origin
	and take $\varepsilon > 0$ such that
	$\bbS_{r}$ intersects $X^{\ast}$ transversely for any $0 < r \leq \varepsilon$.
	We may assume that $M = X^{\ast} \cap \bbS_{r}$ for some $0 < r < \varepsilon$.
	Define
	\begin{equation}
		W_{0} \coloneqq X \cap \bbB_{\varepsilon}, \qquad
		W \coloneqq X \cap \bbB_{r} \Subset W_{0}.
	\end{equation}
	Suppose to the contrary that $M$ admits a pseudo-Einstein structure.
	Then we can take $0 < r^{\prime} < r$ such that
	the canonical bundle $K_{C}$ of
	\begin{equation}
		C \coloneqq X \cap (\bbB_{r} \setminus \overline{\bbB}_{r^{\prime}}) \subset W_{0} \setminus \{0\}
	\end{equation}
	has a flat Hermitian metric.
	Note that $C$ is homotopy equivalent to $M$;
	in particular,
	$C$ is connected,
	and $W \setminus C = X \cap \overline{\bbB}_{r^{\prime}}$ is compact.
	Parallel transport along curves in $C$ defines a holonomy character
	\begin{equation}
		\chi \colon \pi_{1}(C) \longrightarrow U(1).
	\end{equation}
	Since $U(1)$ is abelian,
	this character factors through $H_{1}(C, \bbZ) \cong H_{1}(M, \bbZ)$.
	Choose a positive integer $m$ annihilating $H_{1}(M, \bbZ)$.
	Then $\chi^{m} = 1$,
	so $K_{C}^{\otimes m}$ has a global parallel frame $s$.
	Since the $(0, 1)$-part of the Chern connection is $\delb$,
	the parallel frame $s$ is holomorphic.
	Its dual frame $t$ is also holomorphic,
	and satisfies $t(s) = 1$.
	
	Both $\calE = K_{W_{0}}^{[m]}$ and its dual $\calE^{\vee}$
	are coherent reflexive sheaves.
	By \cref{lem:reflexive-compact-hartogs},
	$s$ and $t$ extend to $W$.
	Their pairing remains equal to one by the identity theorem.
	Extending both the frame and its dual is essential here;
	extending the frame alone would not show that the extended reflexive sheaf is locally free at $o$.
	Consequently,
	the maps
	\begin{equation}
		\calO_{W} \xrightarrow{1 \mapsto s} \calE \stackrel{t}{\longrightarrow} \calO_{W}
	\end{equation}
	have the identity as their composite,
	and $\calE \cong \calO_{W} \oplus \Ker t$.
	Since $\calE$ is a reflexive rank one sheaf,
	$\Ker t$ is a torsion-free sheaf of rank zero,
	so it is equal to zero.
	Thus $K_{W}^{[m]} \cong \calO_{W}$,
	contradicting the assumption that $K_{X}$ is not $\bbQ$-Cartier at $o$.
	When $H_{1}(M, \bbZ)=0$,
	we may take $m = 1$,
	so $K_{X}$ is Cartier at $o$,
	contradicting the second hypothesis.
\end{proof}

\section{A counterexample on an integral homology 3-sphere}
\label{section:counterexample-integral-homology-sphere}

We first apply \cref{prop:finite-holonomy-canonical} in complex dimension two.
For a normal surface germ $(X, o)$,
normality implies the Cohen--Macaulay property,
so the dualizing sheaf agrees with the reflexive canonical sheaf $K_{X}$.
Consequently, the Gorenstein condition is equivalent to invertibility of $K_{X}$ near $o$,
or to $K_{X}$ being Cartier at $o$.

\begin{proposition}
\label{prop:surface-homology-sphere}
	Let $(X, o) \subset (\bbC^{4}, 0)$ be the surface germ defined by the $2 \times 2$ minors of
	\begin{equation}
	\label{eq:surface-homology-sphere-matrix}
		\begin{pmatrix}
			x & y & z \\
			y - 3 t^{4} & z + t^{6} & x^{2} + 6 t^{2} y - 2 t^{6}
		\end{pmatrix}.
	\end{equation}
	This is the normal non-Gorenstein singularity of
	N\'{e}methi--Okuma~\cite{Nemethi-Okuma2017}*{Theorem 4.1 and Example 4.3}.
	Its link $M$ is an integral homology $3$-sphere.
	Although $T^{1, 0} M$ is trivial as a smooth complex line bundle,
	$M$ admits no pseudo-Einstein structure.
\end{proposition}

\begin{proof}
	The non-Gorenstein assertion is the conclusion of
	\cite{Nemethi-Okuma2017}*{Theorem 4.1} for the example constructed in their Example 4.3.
	This implies that $K_{X}$ is not Cartier at $o$.
	For the topology of its link,
	consider a good resolution,
	whose exceptional divisor has smooth components $E_{i}$ with simple normal crossings.
	Its dual graph $\Gamma$ has one vertex for each $E_{i}$,
	labeled by its genus and self-intersection number,
	and one edge for each intersection point of two components.
	Write $I = (E_{i} \cdot E_{j})$ for the intersection matrix.
	In the present example,
	the minimal good resolution graph is the tree of rational curves in
	\cite{Nemethi-Okuma2017}*{Figure 1 and Example 4.3}.
	It is displayed in Figure~\ref{fig:surface-resolution-graph}.
	
	\begin{figure}[htbp]
		\centering
		\begin{tikzpicture}[
			x = 1.3cm, y = 1cm,
			vertex/.style = {circle, fill, inner sep = 1.7pt},
			line width = 0.5pt
		]
			\draw (0, 0) -- (4, 0);
			\draw (1, 0) -- (1, -0.9);
			\draw (3, 0) -- (3, -0.9);
			\node[vertex, label = above:{$-3$}] at (0, 0) {};
			\node[vertex, label = above:{$-1$}] at (1, 0) {};
			\node[vertex, label = above:{$-13$}] at (2, 0) {};
			\node[vertex, label = above:{$-1$}] at (3, 0) {};
			\node[vertex, label = above:{$-3$}] at (4, 0) {};
			\node[vertex, label = below:{$-2$}] at (1, -0.9) {};
			\node[vertex, label = below:{$-2$}] at (3, -0.9) {};
		\end{tikzpicture}
		\caption{The minimal good resolution graph of \cref{prop:surface-homology-sphere}
			in CR dimension $n = 1$, following \cite{Nemethi-Okuma2017}*{Figure 1}.
			Each vertex represents a smooth rational curve, and its weight is the self-intersection number.}
		\label{fig:surface-resolution-graph}
	\end{figure}
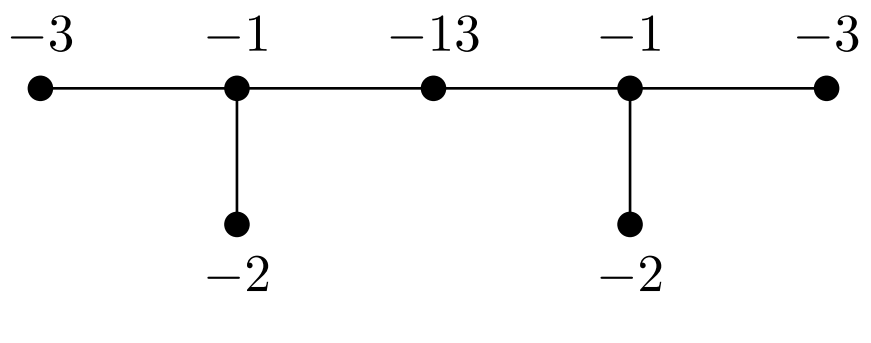
	
	The intersection matrix is
	\begin{equation}
		I =
		\begin{pmatrix}
			-3&1&0&0&0&0&0\\
			1&-1&1&0&0&1&0\\
			0&1&-13&1&0&0&0\\
			0&0&1&-1&1&0&1\\
			0&0&0&1&-3&0&0\\
			0&1&0&0&0&-2&0\\
			0&0&0&1&0&0&-2
		\end{pmatrix}.
	\end{equation}
	For a tree of rational exceptional curves,
	the first homology of the link is presented by
	this intersection matrix \cite{Mumford1961}*{Section I, pp.~10--11}.
	A calculation shows that $\det I = - 1$,
	and so $I \colon \bbZ^{7} \to \bbZ^{7}$ is an isomorphism.
	Consequently,
	\begin{equation}
		H_{1}(M, \bbZ) \cong \operatorname{Coker} I = 0.
	\end{equation}
	Thus $M$ is an integral homology sphere.
	In particular,
	$H^{2}(M, \bbZ) = 0$,
	and hence
	\begin{equation}
		c_{1}(T^{1, 0} M) = 0 \quad \text{in } H^{2}(M, \bbZ).
	\end{equation}
	Since $T^{1, 0} M$ has complex rank one,
	this also implies its triviality as a smooth complex line bundle.
	It follows from $H_{1}(M, \bbZ) = 0$ and \cref{prop:finite-holonomy-canonical} that
	there exist no pseudo-Einstein structures on $M$.
\end{proof}

\section{Line bundle extensions and the first obstruction class}
\label{section:line-bundle-extensions-obstruction}

We explain two constructions in which the differences between local holomorphic choices
define an additive \Cech cocycle.
The first classifies extensions of two line bundles;
the second provides a first obstruction to trivializing a line bundle
on a neighborhood of a hypersurface along which it is already trivial.

\subsection{\Cech classification of line bundle extensions}
\label{subsection:cech-line-bundle-extensions}

Let $B$ be a complex manifold, and let $L$ and $L^{\prime}$ be holomorphic line bundles on $B$.
An \emph{extension of $L$ by $L^{\prime}$} is an exact sequence of holomorphic vector bundles
\begin{equation}
\label{eq:general-line-bundle-extension}
	0 \longrightarrow L^{\prime} \stackrel{\iota}{\longrightarrow} E
	\stackrel{p}{\longrightarrow} L \longrightarrow 0.
\end{equation}
Two extensions are equivalent if there is a holomorphic isomorphism between their middle bundles
that commutes with the inclusions and projections and induces the identity on both $L$ and $L^{\prime}$.
Thus the classification keeps track of these maps as well as the middle bundle.

The sequence splits locally:
a local frame of $L$ has a local holomorphic lift to $E$,
which defines a right inverse of $p$ after shrinking the open set.
Choose an open cover $\{U_{i}\}$ and splittings
\begin{equation}
	s_{i} \colon L|_{U_{i}} \longrightarrow E|_{U_{i}},
	\qquad p \circ s_{i} = \id_{L|_{U_{i}}}.
\end{equation}
On $U_{i} \cap U_{j}$,
the difference $s_{i} - s_{j}$ takes values in $\iota(L^{\prime})$.
Consequently,
there is a unique section
\begin{equation}
	\eta_{i j} \in \Gamma(U_{i} \cap U_{j}, L^{- 1} \otimes L^{\prime})
\end{equation}
satisfying $s_{i} - s_{j} = \iota \circ \eta_{i j}$,
where a section of $L^{- 1} \otimes L^{\prime}$ is viewed as a bundle map from $L$ to $L^{\prime}$.
On a triple overlap,
subtracting the splittings gives
\begin{equation}
	\eta_{i j} + \eta_{j k} = \eta_{i k}.
\end{equation}
This is the additive \Cech cocycle condition.
If the splittings are replaced by $s^{\prime}_{i} = s_{i} + \iota \circ \tau_{i}$,
where $\tau_{i} \in \Gamma(U_{i}, L^{- 1} \otimes L^{\prime})$,
then
\begin{equation}
	\eta^{\prime}_{i j} = \eta_{i j} + \tau_{i} - \tau_{j}.
\end{equation}
Thus changing the local choices adds a coboundary.
Passing to refinements gives a well-defined class
$[\{\eta_{i j}\}] \in H^{1}(B, L^{- 1} \otimes L^{\prime})$.
Equivalent extensions have the same class,
since an equivalence carries local splittings to local splittings.

Conversely,
start with an $L^{- 1} \otimes L^{\prime}$-valued \Cech cocycle $\{\eta_{i j}\}$.
On each $U_{i}$ take a copy $E_{i}$ of $(L^{\prime} \oplus L)|_{U_{i}}$.
On an overlap,
identify a vector in the $i$-th copy with a vector in the $j$-th copy by
\begin{equation}
\label{eq:extension-gluing-from-cocycle}
	(l^{\prime}, l)_{i} \sim (l^{\prime} + \eta_{i j}(l), l)_{j}.
\end{equation}
These are holomorphic bundle isomorphisms,
and the cocycle condition makes the identifications consistent on triple overlaps.
They therefore glue the $E_{i}$ to a rank-two holomorphic bundle $E_{\eta}$.
The local maps $l^{\prime} \mapsto (l^{\prime}, 0)$
and $(l^{\prime}, l) \mapsto l$ also glue,
giving an extension of $L$ by $L^{\prime}$.

If $\eta^{\prime}_{i j} = \eta_{i j} + \tau_{i} - \tau_{j}$,
the local maps
\begin{equation}
	(l^{\prime}, l)_{i} \longmapsto (l^{\prime} - \tau_{i}(l), l)_{i}
\end{equation}
are compatible with the two gluing rules and define an equivalence
$E_{\eta} \cong E_{\eta^{\prime}}$.
The two constructions are inverse: for an extension with splittings $s_{i}$,
the maps $(l^{\prime}, l)_{i} \mapsto \iota(l^{\prime}) + s_{i}(l)$ identify
its reconstruction with the original extension.
We obtain the classification
\begin{equation}
\label{eq:general-extension-classification}
	\Ext^{1}_{\calO_{B}}(L, L^{\prime})
	\cong H^{1}(B, L^{- 1} \otimes L^{\prime}),
\end{equation}
where the left-hand side denotes equivalence classes of extensions with the two end bundles fixed.
The zero class corresponds precisely to a holomorphically split extension:
if $\eta_{i j} = \tau_{i} - \tau_{j}$, the adjusted splittings
$s_{i} - \iota \circ \tau_{i}$ agree and glue to a global splitting.
In particular,
if $H^{1}(B, L^{- 1} \otimes L^{\prime}) = 0$,
every such extension splits.

For comparison with transition functions,
choose frames $l_{i}$ of $L$ and $l_{i}^{\prime}$ of $L^{\prime}$,
and write
\begin{equation}
	l_{i} = \lambda_{i j} l_{j}, \qquad
	l_{i}^{\prime} = \lambda_{i j}^{\prime} l_{j}^{\prime}, \qquad
	\eta_{i j} = c_{i j} l_{j}^{- 1} \otimes l_{j}^{\prime}.
\end{equation}
Let $v_{i} = s_{i}(l_{i})$ and identify $l_{i}^{\prime}$ with $\iota(l_{i}^{\prime})$.
Then
\begin{equation}
\label{eq:general-extension-frame-transition}
	v_{i} 
	= \lambda_{i j} s_{i}(l_{j})
	= \lambda_{i j} (s_{j}(l_{j}) + (s_{i} - s_{j})(l_{j}))
	= \lambda_{i j} (v_{j} + c_{i j} l_{j}^{\prime}).
\end{equation}
The coefficient $c_{i j}$ satisfies
\begin{equation}
	c_{i k} = \frac{\lambda_{j k}^{\prime}}{\lambda_{j k}} c_{i j} + c_{j k}.
\end{equation}
The factor $\lambda_{j k}^{\prime} / \lambda_{j k}$ expresses
the change of frame in $L^{- 1} \otimes L^{\prime}$.

\subsection{The first obstruction to triviality near a hypersurface}
\label{subsection:first-obstruction-triviality}

In the projective bundle construction of \cite{Gongyo-Hirachi-Matsumoto-Takeuchi2026}*{Proposition 2.5},
the obstruction is the class of the line bundle $F|_{D_{1}}$,
where $D_{1}$ is the first infinitesimal neighborhood of the exceptional hypersurface $D$.
Here we construct a \Cech representative of the same obstruction directly on a smooth hypersurface $D$ in a complex manifold $Y$.
We choose local frames that agree on $D$ and take the first normal derivatives of their transition functions.
This construction specializes Koike's first obstruction class
to the case in which the restricted line bundle is holomorphically trivial,
up to the sign determined by our convention for transition functions;
see \cite{Koike2021}*{Section 3.2 and Proposition 3.2(i)}.
No flatness assumption on $N_{D / Y}$ is needed for this first-order construction.

Let $D$ be a smooth connected compact hypersurface in a complex manifold $Y$,
and let $L$ be a holomorphic line bundle on a neighborhood of $D$ with $L|_{D} \cong \calO_{D}$.
Fix a nowhere-vanishing holomorphic section $e$ of $L|_{D}$.
Take a family $\{U_{i}\}$ of Stein open sets in $Y$ such that
$U_{i} \cap D \neq \emptyset$ and $D \subset \bigcup U_{i}$.
Then $\{V_{i} \coloneqq D \cap U_{i}\}$ is a Stein open covering of $D$.
Shrinking $\{U_{i}\}$ if necessary,
we can choose local frames $e_{i}$ of $L$ on each $U_{i}$ such that $e_{i}|_{V_{i}} = e|_{V_{i}}$.
Write
\begin{equation}
	e_{i} = h_{i j} e_{j}, \qquad h_{i j}|_{D} = 1.
\end{equation}
The conormal bundle of $D$ in $Y$ is
\begin{equation}
	N_{D / Y}^{\ast} = \Ker \rbra*{\Omega^{1}_{Y}|_{D} \longrightarrow \Omega^{1}_{D}}.
\end{equation}
Since $h_{i j}$ is constant on $D$,
the ambient differential $d h_{i j}$ annihilates tangent vectors in $T^{1, 0} D$.
We may therefore regard
\begin{equation}
\label{eq:first-normal-cocycle}
	\theta_{i j} = (d h_{i j})|_{D}
	\in \Gamma(V_{i} \cap V_{j}, N_{D / Y}^{\ast})
\end{equation}
as a conormal-valued section.
Differentiating $h_{i j} h_{j k} = h_{i k}$ on $Y$ and then restricting to $D$ gives
\begin{equation}
	\theta_{i j} + \theta_{j k} = \theta_{i k}.
\end{equation}
Thus $\{\theta_{i j}\}$ is an additive \Cech $1$-cocycle.
Passing to refinements if necessary,
we obtain a cohomology class
\begin{equation}
\label{eq:first-order-obstruction-class}
	u_{1}(L) = [\{\theta_{i j}\}] \in H^{1}(D, N_{D / Y}^{\ast}).
\end{equation}

To read this definition in coordinates,
choose local holomorphic coordinates $(z', t_{j})$ in which $D = \{t_{j} = 0\}$.
Locally,
\begin{equation}
	h_{i j} = 1 + b_{i j} t_{j} + O(t_{j}^{2}),
	\qquad
	\theta_{i j} = (b_{i j}|_{D}) \, d t_{j}|_{D}.
\end{equation}
Here $O(t_{j}^{2})$ means an actual holomorphic function divisible by $t_{j}^{2}$.
In particular,
$d t_{j}|_{D}$ is a local frame of the conormal bundle of $D$.

Another family of frames extending the same $e$ has the form
$e^{\prime}_{i} = \lambda_{i} e_{i}$,
where $\lambda_{i}|_{D} = 1$.
Its transition functions satisfy $h^{\prime}_{i j} = \lambda_{i} h_{i j} \lambda_{j}^{- 1}$,
so
\begin{equation}
	\theta^{\prime}_{i j} = \theta_{i j} + (d \lambda_{i})|_{D} - (d \lambda_{j})|_{D}.
\end{equation}
The change is a coboundary,
so the class is independent of these frames.
It is also independent of the cover and of the local coordinates used to compute it.
Any two nowhere-vanishing sections of $L|_{D}$ differ by
a nowhere-vanishing global holomorphic function on $D$.
Since $D$ is connected and compact,
such a function is a non-zero constant;
multiplying every $e_{i}$ by this constant leaves all $h_{i j}$ unchanged.
Thus $u_{1}(L)$ is independent of the initial trivialization of $L|_{D}$ as well.

If $u_{1}(L) = 0$,
then, after refinement,
$\theta_{i j} = \beta_{i} - \beta_{j}$ for local conormal sections $\beta_{i}$.
Choose local holomorphic functions $b_{i}$ vanishing along $D$
with $(d b_{i})|_{D} = - \beta_{i}$.
Replacing $e_{i}$ by $(1 + b_{i})e_{i}$,
after shrinking so that $1 + b_{i} \neq 0$,
makes all first normal derivatives of the transition functions vanish.
This equivalence concerns only the linear terms;
it does not assert that $L$ is trivial on a neighborhood when $u_{1}(L) = 0$.

Multiplication of transition functions and the product rule give
\begin{equation}
\label{eq:first-order-additivity}
	u_{1}(L_{1} \otimes L_{2}) = u_{1}(L_{1}) + u_{1}(L_{2}),
	\qquad
	u_{1}(L^{k}) = k u_{1}(L) \quad (k \in \bbZ),
\end{equation}
whenever the restrictions of the bundles to $D$ are trivial.

\begin{lemma}
\label{lem:first-order-non-torsion}
	Assume that $u_{1}(L) \neq 0$.
	Then $L^{k}|_{\Omega} \not\cong \calO_{\Omega}$ for every $k \in \bbZ \setminus \{0\}$
	and for every open neighborhood $\Omega$ of $D$ contained in the domain of $L$.
\end{lemma}

\begin{proof}
	Suppose to the contrary that $L^{k}|_{\Omega} \cong \calO_{\Omega}$
	for some $k \in \bbZ \setminus \{0\}$
	and for some open neighborhood $\Omega$ of $D$ contained in the domain of $L$.
	Then $L^{k}|_{\Omega}$ has a holomorphic trivialization $\wte$.
	The restriction $e$ of $\wte$ to $D$ is also a global frame on $D$,
	and choosing $e_{j} \coloneqq \wte$ for every $j$ gives $h_{i j} = 1$.
	Hence we have $k u_{1}(L) = u_{1}(L^{k}) = 0$.
	Since $H^{1}(D, N_{D / Y}^{\ast})$ is a complex vector space,
	this means $u_{1}(L) = 0$,
	which is a contradiction.
\end{proof}

\section{Construction of the projective base}
\label{section:projective-base}

We construct a smooth projective base $B$ with an ample line bundle $A$.
The base will become the exceptional hypersurface.
The condition $H^{1}(B, A) \neq 0$ supplies a non-zero extension class $\xi$ for the projective bundle construction.
The vanishing of $H^{1}(B, \calO_{B})$ gives $b_{1}(B) = 0$ by Hodge decomposition.
The substantive existence result is due to Catanese.
There is a connected smooth projective surface $S$ with an ample line bundle $L$ such that
\begin{equation}
\label{eq:catanese-surface}
	K_{S} \cong L^{4}, \qquad
	H^{1}(S, \calO_{S}) = 0, \qquad
	H^{1}(S, L) \neq 0.
\end{equation}
Indeed,
take the Fermat-curve degree to be $7$ and the cohomology exponent to be $1$
in \cite{Catanese2015}*{Theorem 2}.
The degree $7$ satisfies the required lower bound and the coprimality condition $\gcd(7, 30) = 1$;
its canonical exponent is then $7 - 3 = 4$.
The smooth projective surface and the descended line bundle are described in
\cite{Catanese2015}*{Section 4, especially Remark 8}.
In higher dimension,
we take a suitably chosen smooth hypersurface in $S \times \bbP^{n - 1}$.

\begin{proposition}
\label{prop:projective-input}
	For each integer $n \geq 2$,
	there exist a connected smooth complex projective manifold $B$ of complex dimension $n$,
	an ample line bundle $A$ over $B$,
	and an integer $k \geq 2$ such that
	\begin{equation}
	\label{eq:projective-input}
		K_{B} \cong A^{k}, \qquad
		H^{1}(B, \calO_{B}) = 0, \qquad
		H^{1}(B, A) \neq 0.
	\end{equation}
\end{proposition}

\begin{proof}
	In the case of $n = 2$,
	we simply take $(B, A, k) = (S, L, 4)$.
	In what follows,
	we may assume $n \geq 3$,
	and choose $a \geq 2$ so that $L^{a}$ is very ample.
	Consider the product $Z_{n} \coloneqq S \times \bbP^{n - 1}$,
	which is of complex dimension $n + 1$.
	Its canonical bundle $K_{Z_{n}}$ is given by
	\begin{equation}
		K_{Z_{n}} \cong L^{4} \boxtimes \calO(- n).
	\end{equation}
	Since $L^{a}$ and $\calO(n + a + 4)$ are very ample,
	the external product
	\begin{equation}
		L^{a} \boxtimes \calO(n + a + 4)
	\end{equation}
	on $Z_{n}$ is very ample by the Segre embedding.
	Bertini's theorem~\cite{Hartshorne1977}*{Theorem 8.18 in Chapter II} ensures
	the existence of a closed smooth hypersurface $B_{n}$ in $Z_{n}$ such that
	\begin{equation}
		\calO_{Z_{n}}(B_{n}) \cong L^{a} \boxtimes \calO(n + a + 4).
	\end{equation}
	We will verify its connectedness below.
	Set
	\begin{equation}
		H_{n} \coloneqq L \boxtimes \calO(1), \qquad
		A_{n} \coloneqq H_{n}|_{B_{n}}.
	\end{equation}
	The line bundles $H_{n}$ and $A_{n}$ are ample.
	The adjunction formula for $B_{n} \subset Z_{n}$ gives
	\begin{equation}
		K_{B_{n}}
		\cong (K_{Z_{n}} \otimes \calO_{Z_{n}}(B_{n}))|_{B_{n}}
		\cong A_{n}^{a + 4}.
	\end{equation}
	We would like to prove that
	$(B, A, k) = (B_{n}, A_{n}, a + 4)$ satisfies \cref{eq:projective-input}.
	
	Let $i_{n} \colon B_{n} \hookrightarrow Z_{n}$ be the inclusion.
	Restriction to the hypersurface gives an exact sequence
	\begin{equation}
	\label{eq:restriction-to-hypersurface}
		0 \longrightarrow F_{n} \longrightarrow H_{n} \longrightarrow (i_{n})_{\ast} A_{n} \longrightarrow 0,
	\end{equation}
	where 
	\begin{equation}
		F_{n}
		\coloneqq H_{n} \otimes \calO_{Z_{n}}(- B_{n})
		\cong L^{1 - a} \boxtimes \calO(- n - a - 3).
	\end{equation}
	Since $L^{1 - a}$ and $\calO(- n - a - 3)$ are negative,
	we have
	\begin{equation}
		H^{0}(S, L^{1 - a}) = 0, \qquad
		H^{0}(\bbP^{n - 1}, \calO(- n - a - 3)) = 0.
	\end{equation}
	Moreover,
	$H^{p}(\bbP^{n - 1}, \calO(l)) = 0$
	for any $1 \leq p \leq n - 2$ and for any $l \in \bbZ$.
	By the K\"{u}nneth formula,
	\begin{align}
		H^{q}(Z_{n}, F_{n})
		&= \bigoplus_{j = 0}^{q} H^{j}(S, L^{1 - a})
			\otimes_{\bbC} H^{q - j}(\bbP^{n - 1}, \calO(- n - a - 3)) \\
		&= \bigoplus_{j = 1}^{q - 1} H^{j}(S, L^{1 - a})
			\otimes_{\bbC} H^{q - j}(\bbP^{n - 1}, \calO(- n - a - 3)) \\
		&= 0
	\end{align}
	for any $1 \leq q \leq n - 1$.
	In particular,
	$H^{1}(Z_{n}, F_{n}) = 0$ and $H^{2}(Z_{n}, F_{n}) = 0$.
	The long exact sequence associated with \eqref{eq:restriction-to-hypersurface} now gives
	\begin{align}
	\label{eq:nonvanishing-higher-base}
		H^{1}(B_{n}, A_{n})
		\cong H^{1}(Z_{n}, H_{n})
		\cong H^{1}(S, L) \otimes_{\bbC} H^{0}(\bbP^{n - 1}, \calO(1)),
	\end{align}
	whose dimension is $n \dim_{\bbC} H^{1}(S, L) > 0$.
	
	Similarly,
	the restriction of the structure sheaf yields
	\begin{equation}
	\label{eq:restriction-structure-sheaf}
		0 \longrightarrow G_{n} \longrightarrow \calO_{Z_{n}}
		\longrightarrow (i_{n})_{\ast} \calO_{B_{n}} \longrightarrow 0,
	\end{equation}
	where
	\begin{equation}
		G_{n}
		\coloneqq \calO_{Z_{n}}(- B_{n})
		\cong L^{- a} \boxtimes \calO(- n - a - 4).
	\end{equation}
	A similar calculation gives
	\begin{align}
		H^{q}(Z_{n}, G_{n})
		&= \bigoplus_{j = 0}^{q} H^{j}(S, L^{- a})
			\otimes_{\bbC} H^{q - j}(\bbP^{n - 1}, \calO(- n - a - 4)) \\
		&= \bigoplus_{j = 1}^{q - 1} H^{j}(S, L^{- a})
			\otimes_{\bbC} H^{q - j}(\bbP^{n - 1}, \calO(- n - a - 4)) \\
		&= 0
	\end{align}
	for any $1 \leq q \leq n - 1$.
	Consequently,
	the long exact sequence associated with \eqref{eq:restriction-structure-sheaf} now gives
	\begin{equation}
		H^{0}(B_{n}, \calO_{B_{n}}) \cong H^{0}(Z_{n}, \calO_{Z_{n}}) \cong \bbC
	\end{equation}
	and
	\begin{equation}
		H^{1}(B_{n}, \calO_{B_{n}}) \cong H^{1}(Z_{n}, \calO_{Z_{n}}) = 0.
	\end{equation}
	In particular,
	$H^{0}(B_{n}, \bbC) \cong H^{0}(B_{n}, \calO_{B_{n}}) \cong \bbC$ implies that $B_{n}$ is connected.
	Taking $(B, A, k) = (B_{n}, A_{n}, a + 4)$ completes the construction.
\end{proof}

\section{An exceptional hypersurface in a projective bundle}
\label{section:exceptional-hypersurface-projective-bundle}

Let $(B, A, k)$ be as in \cref{prop:projective-input}
and choose $0 \neq \xi \in H^{1}(B, A)$.
This cohomology class defines an extension $V$ of $A^{- 1}$ by $\calO_{B}$:
\begin{equation}
\label{eq:bundle-extension}
	0 \longrightarrow \calO_{B} \stackrel{e_{0}}{\longrightarrow} V
	\longrightarrow A^{- 1} \longrightarrow 0.
\end{equation}
Let $\pi \colon Y \coloneqq \bbP_{B}(V) \to B$ be the projective bundle of $V$;
here we use the quotient convention
\begin{equation}
	Y = \Set{(x, [\ell]) | \ell \in V_{x}^{\ast} \setminus \{0\}}.
\end{equation}
The holomorphic bundle $V$ is algebraic by GAGA~\cite{Hartshorne1977}*{Theorem 2.1 in Appendix B},
and the projective bundle of an algebraic bundle over a projective base is projective~\cite{Hartshorne1977}*{Proposition 7.10(b) in Chapter II}.
Thus $Y$ is a smooth projective manifold of dimension $n + 1$.
This $Y$ has the tautological line bundle
\begin{equation}
	S \coloneqq \Set{(x, [\ell], v) \in \pi^{\ast} V | \ell(v) = 0}
\end{equation}
and the tautological quotient line bundle
\begin{equation}
	Q \coloneqq \Set{(x, \ell, t) | \ell \in V_{x}^{\ast} \setminus \{0\}, \ t \in \bbC} / \sim,
\end{equation}
where $(x, \ell, t) \sim (x, \ell^{\prime}, t^{\prime})$
if there exists $\lambda \in \bbC^{\times}$ such that
$\ell^{\prime} = \lambda \ell$ and $t^{\prime} = \lambda t$.
The vector bundle $\pi^{\ast} V$ is an extension of $Q$ by $S$;
that is,
$\pi^{\ast} V$ fits into the following exact sequence
\begin{equation}
\label{eq:universal-quotient-sequence}
	0 \longrightarrow S \longrightarrow \pi^{\ast} V \stackrel{q}{\longrightarrow} Q \longrightarrow 0.
\end{equation}
The bundle map $q$ is given by
\begin{equation}
	q(v)_{(x, [\ell])} = [(x, \ell, \ell(v))].
\end{equation}
The holomorphic section $e_{0}$ of $V$ defines a closed smooth complex hypersurface
\begin{equation}
	D \coloneqq \Set{(x, [\ell]) \in Y | \ell(e_{0}) = 0}.
\end{equation}
The projection $\pi$ defines a biholomorphism from $D$ to $B$;
in what follows,
we identify $D$ with $B$ via this biholomorphism.
In particular,
$\pi$ is a holomorphic map from $Y$ to $D$ satisfying $\pi|_{D} = \id_{D}$.
Since $B$ is projective,
\begin{equation}
	b_{1}(D)
	= b_{1}(B)
	= 2 \dim_{\bbC} H^{1}(B, \calO_{B})
	= 0.
\end{equation}
Define a holomorphic line bundle $F$ over $Y$ by
\begin{equation}
	F \coloneqq \calO_{Y}(D) \otimes \pi^{\ast} N_{D / Y}^{\ast}.
\end{equation}
Note that
\begin{equation}
	F|_{D} \cong N_{D / Y} \otimes N_{D / Y}^{\ast} \cong \calO_{D}.
\end{equation}

\begin{theorem}
\label{thm:exceptional-hypersurface}
	For each integer $n \geq 2$,
	there exist a smooth complex projective manifold $Y$ of complex dimension $n + 1$,
	a closed connected smooth complex hypersurface $D$ in $Y$ with $b_{1}(D) = 0$,
	a holomorphic map $\pi \colon Y \to D$,
	and an integer $k \geq 2$ such that
	$\pi|_{D} = \id_{D}$,
	$N_{D / Y}$ is negative,
	no non-zero tensor power of $F = \calO_{Y}(D) \otimes \pi^{\ast} N_{D / Y}^{\ast}$
	is holomorphically trivial on any open neighborhood of $D$,
	and
	\begin{equation}
		K_{Y} \otimes \calO_{Y}((k + 1) D) \cong F^{k - 1}.
	\end{equation}
\end{theorem}

\begin{proof}
	We use the above notation.
	We first show that $\calO_{Y}(D) \cong Q$.
	Consider the holomorphic section $q(e_{0})$ of $Q$.
	The definition of $q$ implies
	\begin{equation}
		q(e_{0})_{(x, [\ell])} = [x, \ell, \ell(e_{0})].
	\end{equation}
	The divisor of $q(e_{0})$ is equal to $D$,
	and so $Q$ is holomorphically isomorphic to $\calO_{Y}(D)$.
	
	We next compute the canonical bundle $K_{Y}$.
	The relative tangent bundle $T_{Y / B}$ of $\pi \colon Y \to B$ is given by
	\begin{equation}
		T_{Y / B} \cong \Hom(S, Q) \cong S^{- 1} \otimes Q.
	\end{equation}
	Taking determinants in \eqref{eq:universal-quotient-sequence} yields
	\begin{equation}
		S \otimes Q \cong \pi^{\ast} \det V \cong \pi^{\ast} A^{- 1}
	\end{equation}
	and
	\begin{equation}
	\label{eq:canonical-projective-bundle}
		K_{Y}
		\cong T_{Y / B}^{\ast} \otimes \pi^{\ast} K_{B}
		\cong Q^{- 2} \otimes \pi^{\ast} A^{k - 1}.
	\end{equation}
	
	Choose local frames $a_{i}$ of $A$ on an open cover $\{U_{i}\}$ of $B$
	and take a lift $v_{i}$ of $a_{i}^{- 1}$ to $V$ on $U_{i}$.
	It follows from \eqref{eq:bundle-extension} that
	$(e_{0}, v_{i})$ is a local frame of $V$ on $U_{i}$.
	Set
	\begin{equation}
		U_{i}^{\prime}
		\coloneqq \Set{(x, [\ell]) \in Y | x \in U_{i}, \ \ell(v_{i}) \neq 0}.
	\end{equation}
	This is an open set of $Y$
	and satisfies $D \cap \pi^{- 1}(U_{i}) \subset U_{i}^{\prime}$.
	The fiber coordinate $t_{i}$ on $U_{i}^{\prime}$ is defined by
	\begin{equation}
		t_{i}(x, [\ell])
		= \frac{q(e_{0})_{(x, [\ell])}}{q(v_{i})_{(x, [\ell])}}
		= \frac{\ell(e_{0})}{\ell(v_{i})}.
	\end{equation}
	In this coordinate,
	the hypersurface $D \cap \pi^{- 1}(U_{i})$ is given by $\{t_{i} = 0\}$.
	On the overlap $U_{i} \cap U_{j}$,
	write $a_{i} = g_{i j} a_{j}$ and
	\begin{equation}
	\label{eq:extension-frame-transition}
		v_{i} = g_{i j}^{- 1} (v_{j} + c_{i j} e_{0})
	\end{equation}
	for holomorphic functions $g_{i j}$ and $c_{i j}$;
	$g_{i j}^{- 1}$ appears because $v_{i}$ is a lift of $a_{i}^{- 1}$.
	Hence we have
	\begin{align}
		t_{i}
		= \frac{q(e_{0})}{q(v_{i})}
		= \frac{g_{i j} q(e_{0})}{q(v_{j}) + c_{i j} q(e_{0})}
		= \frac{g_{i j} t_{j}}{1 + c_{i j} t_{j}}.
	\end{align}
	In particular,
	\begin{equation}
		d t_{i}|_{D} = g_{i j} \, d t_{j}|_{D},
	\end{equation}
	which gives an identification of $N_{D / Y}^{\ast}$ and $A$
	by sending $d t_{i}|_{D}$ to $a_{i}$.
	The normal bundle $N_{D / Y}$ is holomorphically isomorphic to $A^{- 1}$,
	which is negative.
	
	Consider the local frames
	\begin{equation}
		f_{i} = q(v_{i}) \otimes \pi^{\ast} a_{i}
	\end{equation}
	of $F \cong Q \otimes \pi^{\ast} A$,
	which satisfy the exact identity
	\begin{equation}
	\label{eq:corrected-transition}
		f_{i} = (1 + c_{i j} t_{j}) f_{j}.
	\end{equation}
	Their restrictions to $D$ glue to a global frame of $F|_{D}$.
	The cocycle \eqref{eq:first-normal-cocycle} is therefore
	\begin{equation}
		\theta_{i j} = c_{i j} \, d t_{j}|_{D} \in \Gamma(U_{i} \cap U_{j}, N_{D / Y}^{\ast}),
	\end{equation}
	which is identified with $c_{i j} a_{j}$.
	Under this identification,
	we have
	\begin{equation}
	\label{eq:corrected-obstruction}
		u_{1}(F) = [\{c_{i j} a_{j}\}] = \xi \neq 0.
	\end{equation}
	It follows from \cref{lem:first-order-non-torsion} that
	no non-zero tensor power of $F$ is holomorphically trivial
	on any open neighborhood of $D$.
	Finally,
	\eqref{eq:canonical-projective-bundle} implies
	\begin{equation}
		K_{Y} \otimes \calO_{Y}((k + 1) D)
		\cong \calO_{Y}((k - 1) D) \otimes \pi^{\ast} N_{D / Y}^{- k + 1}
		\cong F^{k - 1},
	\end{equation}
	which completes the proof.
\end{proof}

In the construction of the counterexample to Lee's conjecture,
we need to compute the first Betti number of the link
and the first Chern class of the canonical bundle.
To this end,
we study the tubular neighborhood of the hypersurface $D$.

\begin{proposition}
\label{prop:top-of-tubular-neighborhood}
	Let $Y$ and $D$ be as in \cref{thm:exceptional-hypersurface}.
	A sufficiently small tubular neighborhood $T$ of $D \subset Y$
	satisfies $H^{1}(T \setminus D, \bbR) = 0$
	and $c_{1}(K_{T \setminus D}) = 0$ in $H^{2}(T \setminus D, \bbZ)$.
\end{proposition}

\begin{proof}
	The punctured neighborhood $T \setminus D$ is homotopy equivalent to
	the unit circle bundle $p \colon P \coloneqq S(N_{D / Y}) \to D$.
	The Euler class of this circle bundle is given by $c_{1}(N_{D / Y})$.
	The Gysin exact sequence yields
	\begin{equation}
		0 \longrightarrow H^{1}(D, \bbR) \stackrel{p^{\ast}}{\longrightarrow} H^{1}(P, \bbR)
		\longrightarrow H^{0}(D, \bbR) \xrightarrow{\smile c_{1}(N_{D / Y})} H^{2}(D, \bbR).
	\end{equation}
	Since $D$ is connected and $N_{D / Y}$ is negative,
	the map $H^{0}(D, \bbR) \to H^{2}(D, \bbR)$ is injective.
	Hence we have
	\begin{equation}
		H^{1}(T \setminus D, \bbR)
		\cong H^{1}(P, \bbR)
		\cong H^{1}(D, \bbR) = 0.
	\end{equation}
	
	Since $D$ is a deformation retract of $T$ and $F|_{D}$ is trivial,
	$c_{1}(F|_{T}) = 0$ in $H^{2}(T, \bbZ)$.
	It follows from $K_{Y} \cong \calO_{Y}(- (k + 1) D) \otimes F^{k - 1}$
	and the fact that $\calO_{Y}(- (k + 1) D)$ is trivial on $Y \setminus D$ that
	\begin{equation}
		c_{1}(K_{T \setminus D}) = c_{1}(K_{Y}|_{T \setminus D}) = 0
	\end{equation}
	in $H^{2}(T \setminus D, \bbZ)$.
\end{proof}

\section{Contraction and counterexamples in higher CR dimensions}
\label{section:contraction-higher-cr-dimensions}

Let $Y$, $D$, $F$, and $k$ be as in \cref{thm:exceptional-hypersurface}.
Since $N_{D / Y}$ is negative,
\cref{thm:grauert-contraction} gives a neighborhood $\Omega$ of $D$
and a proper holomorphic map
\begin{equation}
	q \colon \Omega \to W
\end{equation}
contracting exactly $D$ to a point $o$,
where $W$ is a normal Stein space
and $\Omega \setminus D \to W \setminus \{o\}$ is biholomorphic.
Write $(X, o)$ for the resulting germ.

\begin{proposition}
\label{prop:contracted-canonical-not-q-cartier}
	The normal isolated singularity $(X, o)$ constructed above has a canonical divisor $K_{X}$
	that is not $\bbQ$-Cartier at $o$,
	and its link $M$ satisfies $b_{1}(M) = 0$ and $c_{1}(T^{1, 0} M) = 0$ in $H^{2}(M, \bbZ)$.
	In particular,
	the CR manifold $(M, T^{1, 0} M)$ admits no pseudo-Einstein structure.
\end{proposition}

\begin{proof}
	Fix a local embedding in $\bbC^{N}$ sending $o$ to the origin
	and take $\varepsilon > 0$ such that
	$\bbS_{r}$ intersects $X^{\ast}$ transversely for any $0 < r \leq \varepsilon$.
	We may assume that $M = X^{\ast} \cap \bbS_{r}$ for some $0 < r < \varepsilon$.
	Define
	\begin{equation}
		W_{+} \coloneqq X \cap \bbB_{\varepsilon},
	\end{equation}
	which contains $M$.
	Take a tubular neighborhood $T$ of $D$ so that $T \subset q^{- 1}(W_{+})$.
	The image $q(T)$ is an open neighborhood of $o \in X$
	since
	\begin{equation}
		q(T) = W \setminus q(\Omega \setminus T)
	\end{equation}
	and the right-hand side is open by the properness of $q$.
	Take a sufficiently small $0 < \varepsilon^{\prime} < r$ so that
	\begin{equation}
		W_{-} \coloneqq X \cap \bbB_{\varepsilon^{\prime}} \subset q(T).
	\end{equation}
	Under the identification via $q$,
	we obtain the sequence of inclusions
	\begin{equation}
		U_{-} \coloneqq W_{-} \setminus \{o\} \longrightarrow T \setminus D
		\longrightarrow U_{+} \coloneqq W_{+} \setminus \{o\}.
	\end{equation}
	This induces
	\begin{equation}
		H^{j}(U_{+}, R) \longrightarrow H^{j}(T \setminus D, R)
		\longrightarrow H^{j}(U_{-}, R)
	\end{equation}
	for any coefficient ring $R$.
	On the other hand,
	the composition $H^{j}(U_{+}, R) \to H^{j}(U_{-}, R)$ is an isomorphism
	since the inclusion $U_{-} \hookrightarrow U_{+}$ is a homotopy equivalence.
	If we take $j = 1$ and $R = \bbR$,
	\cref{prop:top-of-tubular-neighborhood} gives $H^{1}(U_{+}, \bbR) = 0$.
	Since $M$ is a deformation retract of $U_{+}$,
	we have $H^{1}(M, \bbR) = 0$ and $b_{1}(M) = 0$.
	If we take $j = 2$ and $R = \bbZ$,
	the image of $c_{1}(K_{U_{+}})$ in $H^{2}(T \setminus D, \bbZ)$ is equal to zero
	by \cref{prop:top-of-tubular-neighborhood}.
	Hence $c_{1}(K_{U_{+}})$ itself is equal to zero in $H^{2}(U_{+}, \bbZ)$.
	Since $c_{1}(T^{1, 0} M) = - c_{1}(K_{U_{+}}|_{M})$,
	we have $c_{1}(T^{1, 0} M) = 0$ in $H^{2}(M, \bbZ)$.
	
	It remains to show that $K_{X}$ is not $\bbQ$-Cartier at $o$.
	Suppose to the contrary that
	$K_{X}^{[m]}$ is invertible at $o$ for some $m > 0$.
	After shrinking $W$ and replacing $\Omega$ by $q^{- 1}(W)$,
	choose a frame $s$ of $K_{W}^{[m]}$.
	On $U = W \setminus \{o\}$,
	this is a nowhere-vanishing holomorphic section of $K_{U}^{\otimes m}$.
	Under the biholomorphism $\Omega \setminus D \cong U$,
	it gives a nowhere-vanishing holomorphic section $\sigma$ of $K_{Y}^{\otimes m}$ on $\Omega \setminus D$.
	
	We would like to extend $\sigma$ meromorphically to all of $\Omega$.
	The smooth complex projective manifold $Y$ admits a non-zero rational section $\tau$ of $K_{Y}^{\otimes m}$:
	over the function field of $Y$,
	this line bundle is a one-dimensional vector space and has a basis.
	The quotient $\sigma / \tau$ is a scalar meromorphic function on $\Omega \setminus D$.
	Transport it to $U$ and extend it meromorphically to $W$ by \cref{thm:stein-hartogs}.
	If $g$ denotes this extension,
	then
	\begin{equation}
		\whxs = (q^{\ast} g) \tau
	\end{equation}
	is the required meromorphic section on $\Omega$.
	
	All zeros and poles of $\whxs$ are supported on $D$.
	Since $D$ is smooth and connected, it is irreducible.
	Thus the divisor of $\whxs$ is $b D$ for some integer $b$, with $b < 0$ in the case of a pole.
	It follows that
	\begin{equation}
		K_{Y}^{\otimes m}|_\Omega \cong \calO_{\Omega}(b D).
	\end{equation}
	Restricting this isomorphism to $D$ gives
	\begin{equation}
		N_{D / Y}^{- m (k + 1)} \cong N_{D / Y}^{b},
	\end{equation}
	because
	\begin{equation}
		K_{Y}|_{D} \cong \calO_{Y}(- (k + 1) D)|_{D} \cong N_{D / Y}^{- (k + 1)}.
	\end{equation}
	Taking first Chern classes and pairing with $c_{1}(N_{D / Y})^{n - 1}$ yields
	\begin{equation}
		(m (k + 1) + b) \int_{D} c_{1}(N_{D / Y})^{n} = 0.
	\end{equation}
	Since $N_{D / Y}$ is negative,
	this forces $b = - m (k + 1)$.
	Consequently,
	\begin{equation}
		F^{m (k - 1)}|_{\Omega}
		\cong \rbra*{K_{Y}^{\otimes m} \otimes \calO_{Y}(m (k + 1) D)}|_{\Omega}
		\cong \calO_{\Omega},
	\end{equation}
	which contradicts \cref{thm:exceptional-hypersurface}.
	The point $o$ is indeed an isolated singularity
	because $K_{X}$ is not invertible at $o$.
	Since $M$ is a compact smooth manifold, $H_{1}(M, \bbZ)$ is finitely generated.
	Its rank is $b_{1}(M) = 0$, so $H_{1}(M, \bbZ)$ is finite.
	By \cref{prop:finite-holonomy-canonical}, $M$ admits no pseudo-Einstein structure.
\end{proof}

\section*{Declaration of AI use}
\label{section:declaration-ai-use}

The author used OpenAI's ChatGPT and Codex in preparing this expository note.
These tools assisted with drafting and revising mathematical explanations and proofs,
organizing the exposition, locating and checking references,
and reviewing the mathematical arguments and English prose.
They also assisted with editing the \LaTeX{} source.
The author reviewed and revised the resulting material
and takes full responsibility for the content of this note,
including its mathematical statements, proofs, and references.


\begin{bibdiv}
\begin{biblist}

\bib{Catanese2015}{article}{
	author={Catanese, Fabrizio},
	title={Subcanonical graded rings which are not Cohen--Macaulay},
	book={
		title={Recent advances in algebraic geometry},
		series={London Mathematical Society Lecture Note Series},
		volume={417},
		publisher={Cambridge University Press},
		address={Cambridge},
	},
	date={2015},
	pages={92--101},
	doi={10.1017/CBO9781107416000.007},
	note={With an appendix by Jonathan Wahl.},
}

\bib{Dong-Wang2026}{misc}{
	author={Dong, Yuxin},
	author={Wang, Haohao},
	title={A counterexample to Lee's conjecture on pseudo-Einstein structures},
	date={2026},
	note={\href{https://arxiv.org/abs/2611.02718}{arXiv:2611.02718}},
}

\bib{Gongyo-Hirachi-Matsumoto-Takeuchi2026}{misc}{
	author={Gongyo, Yoshinori},
	author={Hirachi, Kengo},
	author={Matsumoto, Shuya},
	author={Takeuchi, Yuya},
	title={Pseudo-Einstein Structures on Links of Isolated Singularities},
	date={2026},
	note={\href{https://arxiv.org/abs/2611.04851}{arXiv:2611.04851}},
}

\bib{Grauert1962}{article}{
	author={Grauert, Hans},
	title={\"Uber Modifikationen und exzeptionelle analytische Mengen},
	journal={Math. Ann.},
	volume={146},
	date={1962},
	pages={331--368},
	doi={10.1007/BF01441136},
}

\bib{Grauert-Remmert1984}{book}{
	author={Grauert, Hans},
	author={Remmert, Reinhold},
	title={Coherent analytic sheaves},
	series={Grundlehren der mathematischen Wissenschaften},
	volume={265},
	publisher={Springer-Verlag},
	address={Berlin},
	date={1984},
	doi={10.1007/978-3-642-69582-7},
}

\bib{Grauert-Remmert2004}{book}{
	author={Grauert, Hans},
	author={Remmert, Reinhold},
	title={Theory of Stein spaces},
	series={Classics in Mathematics},
	publisher={Springer-Verlag},
	address={Berlin},
	date={2004},
	note={Translated from the German by Alan Huckleberry; reprint of the 1979 translation},
	doi={10.1007/978-3-642-18921-0},
}

\bib{Hartshorne1977}{book}{
	author={Hartshorne, Robin},
	title={Algebraic geometry},
	series={Graduate Texts in Mathematics},
	volume={52},
	publisher={Springer-Verlag},
	address={New York--Heidelberg},
	date={1977},
	doi={10.1007/978-1-4757-3849-0},
}

\bib{Hirachi-Marugame-Matsumoto2017}{article}{
	author={Hirachi, Kengo},
	author={Marugame, Taiji},
	author={Matsumoto, Yoshihiko},
	title={Variation of total $Q$-prime curvature on CR manifolds},
	journal={Adv. Math.},
	volume={306},
	date={2017},
	pages={1333--1376},
	doi={10.1016/j.aim.2016.11.005},
}

\bib{Kasahara1965}{article}{
	author={Kasahara, Kenkiti},
	title={On Hartogs-Osgood's theorem for Stein spaces},
	journal={J. Math. Soc. Japan},
	volume={17},
	number={3},
	date={1965},
	pages={297--312},
	doi={10.2969/jmsj/01730297},
}

\bib{Koike2021}{article}{
	author={Koike, Takayuki},
	title={Linearization of transition functions of a semi-positive line bundle along a certain submanifold},
	journal={Ann. Inst. Fourier (Grenoble)},
	volume={71},
	number={5},
	date={2021},
	pages={2237--2271},
	doi={10.5802/aif.3439},
}

\bib{Lee1988}{article}{
	author={Lee, John M.},
	title={Pseudo-Einstein structures on CR manifolds},
	journal={Amer. J. Math.},
	volume={110},
	number={1},
	date={1988},
	pages={157--178},
	doi={10.2307/2374543},
}

\bib{Mumford1961}{article}{
	author={Mumford, David},
	title={The topology of normal singularities of an algebraic surface and a criterion for simplicity},
	journal={Publ. Math. Inst. Hautes \'Etudes Sci.},
	volume={9},
	date={1961},
	pages={5--22},
	doi={10.1007/BF02698717},
}

\bib{Nemethi-Okuma2017}{article}{
	author={N\'emethi, Andr\'as},
	author={Okuma, Tomohiro},
	title={Analytic singularities supported by a specific integral homology sphere link},
	journal={Methods Appl. Anal.},
	volume={24},
	number={2},
	date={2017},
	pages={303--320},
	doi={10.4310/MAA.2017.v24.n2.a7},
}

\bib{Seade2019}{article}{
	author={Seade, Jos\'e},
	title={On Milnor's fibration theorem and its offspring after 50 years},
	journal={Bull. Amer. Math. Soc. (N.S.)},
	volume={56},
	number={2},
	date={2019},
	pages={281--348},
	doi={10.1090/bull/1654},
}






\end{biblist}
\end{bibdiv}
\end{document}